\documentclass [a4paper,12pt] {amsart}
\usepackage [latin1]{inputenc}
\usepackage[all]{xy}
\usepackage{amsmath}
\usepackage{amsfonts}
\usepackage{latexsym}  
\usepackage{amssymb}
\usepackage{graphics}
\usepackage{graphicx}
\newtheorem{teorema}{Theorem}[section]
\newtheorem{Lemma}[teorema]{Lemma}
\newtheorem{propos}[teorema]{Proposition}
\newtheorem{corol}[teorema]{Corollary}
\newtheorem{ex}{Example}[section]
\newtheorem{rem}{Remark}[section]
\newtheorem{defin}[teorema]{Definition}
\def\bt{\begin{teorema}}
\def\et{\end{teorema}}
\def\bp{\begin{propos}}
\def\ep{\end{propos}}
\def\bl{\begin{Lemma}}
\def\el{\end{Lemma}}
\def\bc{\begin{corol}}
\def\ec{\end{corol}}
\def\br{\begin{rem}\rm}
\def\er{\end{rem}}
\def\bex{\begin{ex}\rm}
\def\eex{\end{ex}}
\def\bd{\begin{defin}}
\def\ed{\end{defin}}
\def\demo{\par\noindent{\bf Proof.\ }}
\def\enddemo{\ $\Box$\par\vskip.6truecm}
\def\R{{\mathbb R}}   \def\a {\alpha} 
\def\N{{\mathbb N}}     \def\d {\delta} \def\e{\varepsilon}
\def\C{{\mathbb C}}      \def\l{\lambda}

\def\P{{\mathbb P}}
\def\Re{{\sf Re}}
 \def\oli{\overline}

\def\O{\Omega}

\begin {document}
\title[Semi-global extension]{Semi-global extension of maximally complex submanifolds}
\author[Giuseppe Della Sala, Alberto Saracco]{Giuseppe Della Sala,\\ Alberto Saracco}
               \address{G.\ Della Sala, A.\ Saracco: Scuola Normale Superiore, Piazza dei Cavalieri, 7 - I-56126 Pisa, Italy}\email{g.dellasala@sns.it, a.saracco@sns.it}
               \address{G.\ Della Sala: Dipartimento di Matematica, Università di Pisa, Largo Pontecorvo, 5 - I-56127 Pisa, Italy}
               
               \date{\today}

\begin{abstract}Let $A$ be a domain of the boundary of
a strictly pseudoconvex domain $\Omega$ of $\C^n$ and $M$ a smooth,
closed, maximally complex submanifold of  $A$. We find a subdomain
$\widetilde A$ of $\O$, depending only on $\O$ and $A$, and a
complex variety $W\subset \widetilde A$ such that $bW\cap A = M$. Moreover, a generalization to analytic sets of depth at least $4$ is given.
               \end{abstract}
 \maketitle

\section{Introduction}
In the last fifty years, the \emph{boundary problem}, i.e.\ the
problem of characterizing real submanifolds which are boundaries of
\lq\lq something\rq\rq\ analytic, has been widely treaten.

The first result of this kind is due to Wermer \cite{We}: compact
real curves in $\C^n$ are boundaries of complex varieties if and
only if they satisfy a global integral condition, the \emph{moments
condition}. For greater dimension the problem was solved, by Harvey
and Lawson \cite{HL}, proving that an obviously necessary condition
(\emph{maximal complexity}) is also sufficient for compact manifolds
in $\C^n$. Later on, characterizations for closed (non necessarily
compact) submanifolds in $q$-concave open subsets of $\C\P^n$ were
provided by Dolbeault-Henkin and Dihn in \cite{DH, DH2, Di}. A new
approach to the problem in $\C\P^n$ has been recently set forth by
Harvey-Lawson \cite{HL3,HL4, HL5, HL6}.\vspace{0,2cm}

Our goal is to drop the compactness hypothesis. The results in
\cite{DS} deal with the global situation of submanifolds contained
in the boundary of a special class of strongly pseudoconvex
unbounded domains in $\C^n$. In this paper we deal with the boundary
problem for complex analytic varieties in a \lq\lq semi-global\rq\rq\
setting.

More precisely, let $\O\subset\C^n$ be a strongly pseudoconvex open
domain in $\C^n$, and $b\Omega$ its boundary. Let $M$
be a maximally complex ($2m+1$)-dimensional real closed submanifold ($m\geq1$) of some open domain
$A\subset b\O$, and let $K$ be its boundary. We want to find a
domain $\widetilde A$ in $\O$, independent of $M$, and a complex
subvariety $W$ of $\widetilde A$
such that:
\begin{itemize}
    \item[(i)] $b \widetilde A\cap b\Omega=A$;
    \item[(ii)] $bW\cap b\Omega = M$,
\end{itemize}
%which means that $\widetilde A$ is a \lq\lq removable singularity\rq\rq\ for the problem of extension of complex varieties (or of the extension of a maximally complex manifold to a complex variety).

In this paper we show that, if $A\Subset b\Omega$, the problem we
are dealing with has a solution $(\widetilde A,W)$ whose $\widetilde
A$ can be determined in terms of the envelope $\widehat K$ of $K$ with
respect to the algebra of functions holomorphic in a neighbourhood
of $\oli \O$, i.e.

\begin{itemize}
    \item[] \emph{For any maximally complex
$(2m+1)$-dimensional closed real submanifold $M$ of $A$, $m\geq1$, there exists an $(m+1)$-dimensional complex
variety $W$ in $\O \setminus
\widehat K$, with isolated singularities, such that $bW\cap (A\setminus \widehat K) = M\cap
(A\setminus \widehat K) $.}
\end{itemize}
This result echoes that of Lupacciolu on the extension of CR-functions (see \cite[Theorem 2]{Lu})

If $A$ is not relatively compact, this result can be restated in
terms of \lq\lq principal divisors hull\rq\rq, leading to a global
result for unbounded strictly pseudoconvex domains, different from
the results in \cite{DS}. Indeed, this method of proof allows us to
drop the Lupacciolu hypothesis in \cite{DS} and extend the maximally
complex submanifold to a domain, which can anyhow not be the whole of
$\O$. If the Lupacciolu hypothesis holds, then the domain of
extension is in fact all of $\Omega$. So this result is actually a
generalization of the one in \cite{DS}.

The crucial question of the maximality of the domain $\widetilde A$
we construct is not answered; in some simple cases the domain is
indeed maximal (see Example \ref{ex}).

In the last section, by the same methods, the extension result is
proved for analytic sets (see Theorem \ref{gener}).

It worths noticing that in \cite{ST2} related results are obtained
via a bump Lemma and cohomological methods. That approach may be
generalized to complex spaces.\vspace{0.2cm}

We wish to thank Giuseppe Tomassini for suggesting us the problem in the first place and for useful discussions.

\section{Definitions and notations}
In all the paper we will always consider, unless otherwise stated, $\C^n$ with coordinates
$z_1=x_1+iy_1,\ldots,z_n=x_n+iy_n$,
$x_1,\ldots,x_n,y_1,\ldots,y_n\in\R$.\vspace{0,2cm}

A smooth real ($2m+1$)-dimensional submanifold $M$ of $\C^n$ is said
to be a \emph{CR manifold} if its complex tangent $H_pM$ has
constant dimension at each point $p$. If $m>0$ and $\dim_\C H_pM
=m$, i.e.\ it is the maximal possible, $M$ is said to be
\emph{maximally complex}. Observe that a smooth hypersurface of
$\C^n$ is always maximally complex.

If $m=0$ and $M=\gamma$ is a compact curve, we say that $\gamma$ satisfies the \emph{moments condition} if
$$\int_\gamma \omega\ =\ 0,$$
for any holomorphic $(1,0)$-form $\omega$.

It is easy to observe that the (smooth) boundary of a complex
variety of $\C^n$ of dimension $m+1$ is maximally complex if $m>0$
(respectively satisfies the moments condition if
$m=0$).\vspace{0,2cm}

%Let $\rho:\C^n\to\R$ be a function of class $C^2$. We denote by $\mathcal L_p(\rho)$ the Levi-form of $\rho$ at $p\in\C^n$
%$$\mathcal L_p(\rho)[\xi,\eta]\ =\ \sum_{i,j}\frac{\partial^2\rho}{\partial z_i\partial\oli z_j}\xi_i\oli\eta_j,$$
%where $\xi,\eta\in T^{(1,0)}_p\C^n$. $\rho$ is said to be \emph{strictly plurisubharmonic} if $\mathcal L_p(\rho)$ is strictly positive definite.

A domain $\O\subseteq\C^n$ is called \emph{strongly
pseudoconvex} if there is a neighbourhood $U$ of its boundary such that
$$\O\cap U=\{z\in U:\ \rho(z)<0\},$$
where $\rho$ is strictly plurisubharmonic in $U$.

\section{Main result}
Let $\O\subset\C^n$ be a strongly pseudoconvex open
domain in $\C^n$. Let
$A$ be a subdomain of $b\O$, and $K=b A$. For any Stein neighborhood $\O_\alpha$ of $\O$
we set $\widehat K_\alpha$ to be the hull of $K$ with respect to
the algebra of holomorphic functions of $\O_\alpha$, i.e.\
$$\widehat K_\alpha = \left\{x\in \O_\alpha\ :\ |f(x)|\leq \| f\|_K\
\forall f\in \mathcal O(\O_\alpha)\right\}.$$ We define $\widehat
K$ as the intersection of the $\widehat K_\alpha$ when $\O_\alpha$
varies through the family of all  Stein neighborhoods of $\O$. Observe that,
since $\O$ is strongly pseudoconvex (and thus admits a fundamental
system of Stein neighborhoods, see \cite{T83}), $\widehat K$ coincides with the
hull of $K$ with respect to the algebra of the functions which are
holomorphic in some neighborhood of $\oli\O$. We claim that the
following result holds:
\bt \label{main} For any maximally complex
$(2m+1)$-dimensional closed real submanifold $M$ of $A$, $m\geq1$, there exists an $(m+1)$-dimensional complex
variety $W$ in $\O \setminus
\widehat K$, with isolated singularities, such that $bW\cap (A\setminus \widehat K) = M\cap
(A\setminus \widehat K) $.

%In particular, if $K$ is holomorphically convex, i.e.\ $\widehat K=K$, then the extension is in $\O$ and $bW\cap A=M$.
\et 

Following the same strategy as in \cite{DS} we first have a
semi-local extension result (see Lemma \ref{strip} below). In order to
\lq\lq\ globalize" the extension the main differences with respect to
\cite{DS} are due to the fact that we have to cut $\O$ with level-sets
of holomorphic functions instead of hyperplanes. This creates some
additional difficulties: first of all it is no longer possible to
use the parameter which defines the level-sets as a coordinate;
secondly the intersections between tubular domains (see Lemmas
\ref{37}, \ref{Wg} and \ref{agree}) may not be connected.

With the same proof as in \cite{DS} we have
\bl
\label{strip} There exist a tubular neighborhood $I$ of $A$ in
$\O$ and an $(m+1)$-dimensional complex submanifold with boundary
$W_I\subset\oli\O\cap I$ such that $S\cap bW_I=M$ \el

Now, the hypothesis on the hull of $K$ allows us to prove the following

\bl \label{f=0} Let $z^0\in \O\setminus \widehat K$. Then there
exist an open Stein neighborhood $\O_\alpha \supset \O$ and
$f\in\mathcal O(\O_\alpha)$ such that
\begin{itemize}
  \item [1)]$f(z^0)=0$;
  \item [2)]$\{f=0\}$ is a regular complex hypersurface of $\O_\alpha\setminus\widehat K$;
  \item [3)]$\{f=0\}$ intersects $M$ transversally in a compact manifold.
\end{itemize}\el

\br If $f$ is such a function for $z^0$, for any
point $z'$ sufficiently near to $z^0$, $f(z)-f(z')$ satisfies conditions 1), 2) and 3) for $z'$.
\er

\demo By definition of $\widehat K$, since $z^0\in \O\setminus
\widehat K$ there is a Stein neighborhood $\O_\alpha$ such that
$z^0\not\in \widehat K_\alpha$. So we can find a holomorphic function $g$ in $\O_\alpha$ such that $g(z^0)=1$ and $\|g\|_K < 1$;
$h(z)=g(z) - 1$ is a holomorphic function whose zero set does not
intersect $\widehat K$. Since regular level sets are dense, by
choosing a suitable small vector $v$ and redefining $h$ as
$h(z+v)-h(z^0 + v)$ we can safely assume that $h$ satisfies
both 1) and 2).

We remark that $\{h=0\}\cap b\O\Subset A$ by Alexander's Theorem (see \cite[Theorem 3]{A}), and this shows compactness. Then, we may suppose that $M$ is not contained in $\{z_1 =
z_1^0\}$ and, for $\e$ small enough, we consider the function
$f(z)=h(z) + \e (z_1 - z_1^0)$. It's not difficult to see (by
applying Sard's Lemma) that 3) holds for generic $\e$. \enddemo

Now, we divide the proof of Theorem \ref{main} in two cases: $m\geq2$ and $m=1$. This is
due to the fact that in the latter case proving that we can apply
Harvey-Lawson to $\{f=0\}\cap M$ is not automatic.

\subsection{Dimension of $M$ greater than or equal to $5$: $m\geq2$}\label{sez31}

For any $z^0\in\Omega\setminus\widehat K$, Lemma \ref{f=0} provides a holomorphic function such
that the level-set $f_0=\{f=0\}$ contains $z^0$ and intersects $M$ transversally in
a compact manifold $M_0$. The intersection is again maximally
complex (it is the intersection of a complex manifold and a maximally complex manifold, see \cite{HL}), so we can apply Harvey-Lawson Theorem to obtain
a holomorphic chain $W_0$ such that $bW_0 = M_0$. For
$\tau$ in a small neighborhood $U$ of $0$ in $\C$, the hypersurface $f_\tau =\{f-\tau=0\}$ intersects $M$ transversally along a compact submanifold $M_\tau$ which, again by
Harvey-Lawson Theorem, bounds a holomorphic chain $W_\tau$. Observe that since
$M_\tau\subset f_\tau$, $W_\tau\subset f_\tau$.

We claim the following proposition holds:
\bp
\label{analytic} The union $W_U = \bigcup_{\tau\in U} W_\tau$ is a
complex variety contained in the open set $\widetilde
U=\bigcup_{\tau\in U} f_\tau$.
\ep
We need some intermediate results.  Let us consider a generic
projection $\pi:\widetilde{U}\to\C^{m}$ and set $\C^n=\C^{m+1}\times\C^{n-m-1}$,  with holomorphic coordinates
$(w',w)$, $w'\in\C^{m+1}$, $w=(w_1,\ldots,w_{n-m-1})\in\C^{n-m-1}$. Let
$V_\tau= \C^{m+1} \setminus \pi(M_\tau)$.

For $\tau\in U$,
$w'\in\C^{m+1}\setminus\pi(M_\tau)$ and $\alpha\in\N^{n-m-1}$, we
define $$ I^\alpha(w',\tau)\ \doteqdot\ \int_{(\eta',\eta)\in
M_\tau}\eta^\alpha\omega_{BM}(\eta'-w'),$$
where $\omega_{BM}$ is the
Bochner-Martinelli kernel.

In \cite{DS} the following was proved

\bl \label{Zaitsev}
Let $F(w',\tau)$ be the multiple-valued
function which represents $W_\tau$ on $\C^{m+1}\setminus\pi(M_\tau)$
and denote by $P^\alpha(F(w',\tau))$ the sum of the
$\alpha^{th}$ powers of the values of $F(w',\tau)$. Then
$$
P^\alpha(F(w',\tau)) = I^\alpha(w',\tau).
$$
 In particular,
the cardinality $P^0(F(w',\tau))$ of $F(w',\tau)$ is finite.
\el
\br
Lemma \ref{Zaitsev} implies, in particular, that the functions
$P^\alpha(F(w',\tau))$ are continuous in $\tau$. Indeed, they are
represented as integrals of a fixed form over a submanifold $M_\tau$
which varies continuously with the parameter $\tau$.
\er
\bl\label{Lemma36} $P^\alpha(F(w',\tau))$ is holomorphic in the variable $\tau\in
U\subset\C$, for each $\alpha\in\N^{n-m-1}$. \el
\demo
Let us fix a point $\left(w',\underline{\tau}\right)$ such that $w'\notin
M_{\underline{\tau}}$ (this condition remains true for $\tau\in
B_\epsilon(\underline{\tau})$). Consider as domain of $P^\alpha(F)$
the set $\left\{w'\right\}\times B_\epsilon(\underline{\tau})$. In
view of Morera's Theorem, we need to prove that for any simple
curve $\gamma\subset B_\epsilon(\underline{\tau})$, $$ \int_\gamma
P^\alpha(F(w',\tau))d\tau\ =\ 0. $$ Let $\Gamma\subset
B_\epsilon(\underline{\tau})$ be an open set such that
$b\Gamma=\gamma$. By $\gamma\ast M_\tau$ ($\Gamma\ast M_\tau$) we mean
the union of $M_\tau$ along $\gamma$ (along $\Gamma$). Note that
these sets are submanifolds of $\C\times\C^{n}$. The projection $\pi:\Gamma\ast M_\tau\to\C^n$ on the second factor is injective and $\pi(\Gamma\ast M_\tau)$ is an open
subset of $M$ bounded by $\pi(b\Gamma\ast M_\tau)=\pi(\gamma\ast M_\tau)$. By
Lemma~\ref{Zaitsev} and Stokes Theorem
\begin{eqnarray}
    \nonumber\int_\gamma P^\alpha(F(w',\tau))d\tau\ &=& \int_\gamma I^\alpha(w',\tau)d\tau\ =\\
    \nonumber&=&\ \int_\gamma\left(\int_{(\eta',\eta)\in M_\tau} \eta^\alpha\omega_{BM}(\eta'-w')\right)d\tau\ =\\
    \nonumber&=&\ \iint_{\gamma\ast M_\tau}\eta^\alpha\omega_{BM}(\eta'-w')\wedge d\tau\ =\\
    \nonumber&=&\ \iint_{\Gamma\ast M_\tau}d\left(\eta^\alpha\omega_{BM}(\eta'-w')\wedge d\tau\right)\ =\\
    \nonumber&=&\ \iint_{\Gamma\ast M_\tau}d \eta^\alpha\wedge\omega_{BM}(\eta'-w')\wedge d\tau\ =\\\nonumber&=&\ \iint_{\pi(\Gamma\ast M_\tau)}d \eta^\alpha\wedge\omega_{BM}(\eta'-w')\wedge \pi_\ast d\tau\ =\\
    \nonumber&=& 0.
\end{eqnarray}
The last equality follows from the fact that in $d\eta^\alpha$ appear only holomorphic differentials, $\eta^\alpha$ being holomorphic. But since all the holomorphic differentials supported by $\pi(\Gamma\ast M_\tau)\subset M$ already appear in $\omega_{BM}(\eta'-w')\wedge \pi_\ast d\tau$ (due to the fact that $M$ is maximally complex and contains only $m+1$ holomorphic differentials) the integral is zero.
\enddemo
\noindent
\textbf{Proof of Proposition \ref{analytic}.} From \cite{HL2} it follows  that each
$W_\tau$ has isolated singularities\footnote{There could be singularities coming up from intersections of the solutions relative to different connected components of $M_\tau$. These singularities are analytic sets and therefore should intersect the boundary. This cannot happen and so also these singularities are isolated.}. So, let us fix a regular point
$(w'_0,w_0)\in f_{\tau_0}\subset \widetilde U$. In a neighborhood of
this point $W=W_U$ is a manifold, since the construction depends continuously on the initial data. We want to show that $W$ is indeed
analytic in $\widetilde U$.

 Let us fix $j\in\{1,\ldots,n-m-1\}$ and consider multiindexes $\alpha$ of the form
$(0,\ldots,0,\alpha_j,0,\ldots,0)$; let $P_j^\alpha$ be the
corresponding $P^\alpha(F(w',\tau))$. Observe that for any $j$ we can
consider a finite number of $P_j^\alpha$ (it suffices to use
$h=P_j^0(F(w',\tau))$ of them; not that $h$ is independent of $j$). By a linear combination  of the $P_j^\alpha$ with
rational coefficients, we obtain the
elementary symmetric functions
$S_j^0(w',k),\ldots,S_j^{h}(w',\tau)$ in such a way that for any
point $(w',w)\in W$ there exists $\tau\in U$ such that $(w',w)\in W_\tau$; thus, defining
$$Q_j(w',w,\tau)
= S_j^{h}(w',\tau) + S_j^{h-1}(w',\tau)w_j + \cdots +
S_j^0(w',\tau)w_j^{h} = 0,$$ 
we have, in other words,
$$W\subset V =
\bigcup_{\tau\in U} \bigcap_{j=1}^{n-m-1}\{Q_j(w',w,\tau) = 0\}.$$
Define $\widetilde V\subset \C^n(w',w)\times \C(\tau)$ as
$$
\widetilde V = \bigcap_{j=1}^{n-m-1}\{Q_j(w',w,\tau) = 0\}
$$
and
$$\widetilde W = W_\tau\ast U\subset\widetilde V.$$
Observe that, since the functions $S_j^\alpha$ are holomorphic, $\widetilde
V$ is a complex subvariety of $\C^n \times U$. Since $\widetilde V$ and $\widetilde W$ have the same dimension, in a neighborhood of $(w'_0,w_0,\tau)$ $\widetilde W$
is an open subset of the regular part of $\widetilde V$, thus a complex submanifold. We denote by $\emph{Reg}\ (\widetilde W)$ the set of points $z\in \widetilde W$ such that $\widetilde W\cap{{\mathcal U}}$ is a complex submanifold in a neighborhood ${\mathcal U}$ of $z$. It is easily seen that $\emph{Reg}\ (\widetilde W)$ is an open and closed subset of $\emph{Reg}\ (\widetilde V)$, so a connected component. Observing that the closure of a connected component of the regular part of a complex variety is a complex variety we obtain the that $\widetilde W$ is a complex variety, $\widetilde W$ being the closure of $\emph{Reg}\ (\widetilde W)$ in $\widetilde V$.

Finally, since the projection $\pi:\widetilde W\to W$ is a homeomorphism and so proper, it follows that $W$ is a complex subvariety as well.
\enddemo

Now we prove that the varieties $\widetilde W_U$ that we have found --- which
are defined in the open subsets of type $\widetilde U$ (see
Proposition \ref{analytic}) --- patch together in such a way to define a complex variety on the whole of
$\O\setminus\widehat K$.

\bl\label{37} Let $\widetilde U_f$ and $\widetilde U_g$ be two open subsets as in Proposition \ref{analytic} and let $W_f$ and $W_g$ be the corresponding varieties. Let $z^1\in \widetilde U_f \cap \widetilde U_g$. Then $W_f$ and $W_g$ coincide in a neighborhood of $z^1$. \el
\demo Let $\l=f(z^1)$ and $\tau=g(z^1)$ and consider
$$L(\l',\tau') = \{f=\l'\}\cap\{g=\tau'\}\subset \O$$
for $(\l',\tau')$ in a neighborhood of $(\l, \tau)$. Note that for almost every $(\l',\tau')$ $L(\l',\tau')$ is a complex submanifold of codimension $2$ of $\widetilde U_f\cap \widetilde U_g$. Moreover, $W_f\cap L(\l',\tau')$ and $W_g\cap L(\l',\tau')$ are both solutions of the Harvey-Lawson problem for $M\cap L(\l',\tau')$, consequently they must coincide. Since the complex subvarieties $L(\l',\tau')$ which are regular form a dense subset, $W_f$ and $W_g$ coincide on the connected component of $\widetilde U_f\cap\widetilde U_g$ containing $z^1$. \enddemo

\br
The above proof does not work in the case $m=1$ since $M\cap L(\l',\tau')$ is generically empty.
\er

In order to end the proof of Theorem \ref{main}, we have to show that the set $S$ of the singular points of $W$ is a discrete subset of $\O\setminus \widehat K$. Let $z^1\in \O\setminus \widehat K$, and choose a function $h$, holomorphic in a neighborhood of $\O$ such that $h(z^1)=1$ and $K\subset \{|h|\leq \frac{1}{2}\}$ and consider $f=h-\frac{3}{4}$. Observe that $z^1\in \{\Re  f>0\}$ and $K\subset \{\Re  f<0\}$. Choose a defining function $\varphi$ for $b\O$, strongly psh in a neighborhood of $\O$ and let us consider the family
$$\left(\phi_\lambda=\lambda(\varphi)+(1-\lambda)\Re f\right)_{\lambda\in[0,1]}$$
of strongly plurisubharmonic functions. For $\lambda$ near to $1$, $\left\{\phi_\lambda=0\right\}$ does not intersect the singular locus. Let $\oli \lambda$ be the biggest value of $\lambda$ for which $\{\phi_\lambda=0\}\cap  S \neq\emptyset$. Then the analytic set $S$ touches the boundary of the Stein domain $$\left\{\phi_{\oli \lambda}<0\right\}\cap \O \subset \O.$$ So $\{\phi_{\oli \lambda}=0\}\cap S$ is a set of isolated points in $S$. By repeating the same argument, we conclude that $S$ is made up by isolated points.

\subsection{Dimension of $M$ equal to $3$: $m=1$.}

%Now we have to prove that the varieties that we have found - which are defined in open subsets of the type of $\widetilde U$ (see Proposition \ref{analytic}) - agree on the intersections of two such subsets, thus defining a variety on the whole $\O\setminus\widehat K$. We divide the proof of the previous statement in two Lemmas:

%In this subsection we use the notation introduced in Proposition~\ref{P}.

The first goal is to show that when we slice transversally $M$ with complex hypersurfaces, we obtain $1$-dimensional real submanifolds which satisfy the moments condition.

Again, we fix our attention to a neighborhood of the form
$$\widetilde{U}=\bigcup_{\tau\in U} g_\tau.$$
Let us choose an arbitrary
holomorphic $(1,0)$-form $\omega$ in $\C^n$.

\bl\label{omeol}
The function
$$\Phi_\omega(\tau)\ =\ \int_{M_\tau}\omega$$
is holomorphic in $U$.
\el
\demo
Using again Morera's Theorem, we need to prove that for any simple curve $\gamma\subset U$, $\gamma=b\Gamma$,
$$
\int_\gamma \Phi_\omega(\tau)d\tau\ =\ 0.
$$
By Stokes Theorem, we have
\begin{eqnarray}
    \nonumber\int_\gamma \Phi_\omega(\tau )d\tau \ &=& \int_\gamma\left(\int_{M_k}\omega\right)d\tau \ =\\
    \nonumber &=& \iint_{\gamma\ast M_\tau }\omega\wedge d\tau \ =\\
    \nonumber &=& \iint_{\Gamma\ast M_\tau }d(\omega\wedge d\tau )\ =\\
    \nonumber &=& \iint_{\Gamma\ast M_\tau }\partial\omega\wedge d\tau \ =\\
    \nonumber &=& \iint_{\pi(\Gamma\ast M_\tau )}\partial\omega\wedge \pi_\ast d\tau \ =\\
    \nonumber &=& \ 0.
\end{eqnarray}
The last equality is due to the fact that $\pi(\Gamma\ast M_\tau )\subset M$ is maximally complex and thus supports only $(2,1)$ and $(1,2)$-forms, while $\partial\omega\wedge \pi_\ast d\tau $ is a $(3,0)$-form.
\enddemo
\bl\label{Wg}
Let $g$ be a holomorphic
function on a neighborhood of $\O$, and suppose $\{|g|>1\}\cap
\widehat K=\emptyset$. Then there exists a variety $W_g$ on
$\O\cap\{|g|>1\}$ such that $b W_g \cap b\O=M\cap\{|g|>1\}$.\el
\bl\label{agree} Given two functions $g_1$ and $g_2$ as above,
then $W_{g_1}$ and $W_{g_2}$ agree on
$\{|g_1|>1\}\cap\{|g_2|>1\}$.\el

\noindent{\bf Proof of Lemma \ref{Wg}.}
We are going to use several times open subsets of the type $\widetilde{U}$ as in Proposition \ref{analytic}, so we need to fix some notations. Given  an open subset $U\subset \C$, define $\widetilde{U}$ by
$$\widetilde{U}=\bigcup_{\tau \in U}\{f=\tau \}.$$
From now on we use open subsets of the form $U=B(\oli \tau , \d)$, where $B(\oli \tau , \d)$ is the disc centered at $\oli \tau $ of radius $\d$. We say that $\{f=\oli \tau \}$ is the \textit{core} of $\widetilde{U}$ and $\d$ is its \textit{amplitude}.

For a fixed $d>1$ consider the compact
set $H_d= \oli \O \cap \{|g|\geq d\}$; we show that $W_g$ is well
defined on $H_d$. Let us fix also a compact set $C\subset\O$ such that $W_I$ (see Lemma \ref{strip}) is a closed submanifold in $H_d\setminus C$.

Consider all the open subsets $V_\a=\widetilde{U}_\a\cap \O$, constructed using only the function $f=g-1$ up to addition of the function
$\e (z_j - z_j^0)$ (see Lemma \ref{f=0}). If we
do not allow $\e$ to be greater than a fixed $\oli\e>0$, then by a
standard argument of semicontinuity and compactness we may suppose that
the amplitude of each $\widetilde U$ is greater
than a positive $\d$.

We claim that it is possible to find a countable covering of $H_d$ made by a countable sequence $V_i$ of those $V_\a$ in such a way to have
\begin{enumerate}
    \item $V_0\subset H_d\setminus C$;
    \item if
    $$B_l\ =\ \bigcup_{i=1}^l V_i$$
    then $V_{l+1}\cap B_l\cap \O \neq \emptyset$.
\end{enumerate}

The only thing we have to prove is the existence of $V_0$, since the second statement follows by a standard compactness argument.

Set $L=\max_{H_d} \Re g$. Since $\Re g$ is a non constant pluriharmonic function, $\{\Re g=L\}$ is a compact subset of $b\O \cap H_d$. Then we can choose $\eta>0$ such that $\{\Re g=L-\eta\}\cap \O$ is contained in $H_d\setminus C$, and this allows to define $V_0$.

Let $\widetilde U_1$ and $\widetilde U_2$ be two such open sets
and $z^0\in \widetilde U_1 \cap \widetilde U_2$. We can
suppose that the cores of $\widetilde U_1$ and
$\widetilde U_2$ contain $z^0$. They are of the form $$f+\e_1(z_j -
z_j^0)=\tau(\e_1)\  {\rm and}\  f+\e_2(z_j - z_j^0)=\tau(\e_2).$$ For
$\e\in(\e_1,\e_2)$, we consider the open sets $\widetilde U_\e$
whose core, passing by $z^0$, is $f+\e(z_j - z_j^0)=\tau(\e)$.
We must show that the set
$$\Lambda=\left\{\e\in(\e_1,\e_2): \exists W_\e \ s.t. \  W_\e\cap (\widetilde U_1 \cap \widetilde U_\e) = W_1 \cap (\widetilde U_1 \cap\widetilde U_\e) \right\}$$ is open and closed, where $W_\e$ is a variety in $\widetilde U_\e$.

$\Lambda$ is open. Indeed, if $\e\in \Lambda$, then for $\e'$ in a
neighborhood of $\e$ the core of $\widetilde U_{\e'}$ is contained in $\widetilde U_\e$ and so its intersection with $M$ is maximally complex. Because of Lemma \ref{omeol} the condition holds also for all the level sets in $\widetilde U_{\e'}$ and then we can apply again the Harvey-Lawson Theorem \cite{HL} and the arguments of Proposition \ref{analytic} in order to obtain $W_{\e'}$. Moreover, there is a connected component of $U_\e\cap
U_{\e'}$ which contains $z^0$ and touches the boundary of $\O$,
where the $W_\e$ and $W_{\e'}$ both coincide with $W_I$  (see
Lemma \ref{strip}). By virtue of the analytic continuation principle, they must coincide in the whole
connected component.

$\Lambda$ is closed. Indeed, since each $\widetilde U$ has an amplitude
of at least $\d$, we again have that, for $\oli \e \in\oli \Lambda$, the intersection of $\widetilde U_{\oli
\e}$ and $\widetilde U_\e$ must include (for $\e\in \Lambda$, $|\e - \oli
\e|$ sufficiently small) a connected component containing $z^0$ and touching
the boundary. We then conclude as in the previous case.
\enddemo

\noindent{\bf Proof of Lemma \ref{agree}.} Let us consider the connected components of $W_{g_1}\cap\left\{|g_2|>1\right\}$. For each connected component $W_1$ two cases are possible:
\begin{enumerate}
    \item $W_1$ touches the boundary of $\O$: $W_1\cap b\O\neq\emptyset$;
    \item the boundary of $W_1$ is inside $\O$: $$bW_1\Subset \left\{|g_1|=1\right\}\cup\left\{|g_2|=1\right\}\subset\O$$
\end{enumerate}
In the former, the result easily follows in view of the analytic continuation principle (remember that on a strip near the boundary $W_{g_1}$ and $W_{g_2}$ coincide).

The latter is actually impossible. Indeed, suppose by contradiction that the component $W_1$ satisfies (2). Restrict $g_1$ and $g_2$ to $W_1$ and choose $t>1$ such that
$$W_t\ \doteqdot\ \left\{|g_i|>t,\ i=1,2\right\}\ \Subset\ W_1.$$
The boundary $b W_t$ of $W_t$ consists of points where either $|g_1|=t$ or $|g_2|=t$. Choose a point $z_0$ of the boundary where $|g_1|=t$ and $|g_2|>t$, then $|g_2|$ is a plurisubharmonic function on the analytic set $$A=\left\{g_1=g_1(z_0)\right\}\cap\left\{|g_2|\geq t\right\}.$$ Since $W_t\Subset W_1$, the boundary of the connected component of $A$ through $z_0$ is contained in $\left\{|g_2|=t\right\}$. This is a contradiction, because of the maximum principle for plurisubharmonic functions.
\enddemo

\section{Some remarks}

\subsection{Maximality of the solution}
As stated above, we have not a complete answer to the problem of the
maximality of $\widetilde A$. Nevertheless, here is a simple example
where the constructed domain is actually maximal.

\bex\label{ex} Let $\O\subset \C^n$ be a strongly convex domain with smooth
boundary, $0\in \O$, and let $h$ be a pluriharmonic function defined
in a neighborhood $U$ of $\oli \O$ such that $h(0)=0$ and
$h(z)=h(z_1,\ldots,z_{n-1},0)$ (i.e.\ $h$ does not depend on $z_n$).
Pose
$$H=\{z\in U:\  h(z)=0\}$$
and let
$$A=b\O\cap \{z\in U:\ h(z)>0\}.$$
Then
$$\widetilde A= \O\cap \{z\in U:\ h(z)>0\}.$$
In order to show that $\widetilde A$ is maximal for our problem, it
suffices to find, for any $z\in H\cap \O$, a complex manifold
$W_z\subset \widetilde A$ such that $M_z=\oli W_z \cap A$ is smooth
and $W_z$ cannot be extended through any neighborhood of $z$. We may
suppose $z=0$.

So, let $f\in \mathcal O(\oli \O)$ be such that $\Re f=h$, $f(0)=0$.
We define
$$W_0=\{z\in \widetilde A:\ z_n=e^{\frac{1}{f(z)}}\};$$
$W_0$ extends as a closed submanifold of $U\setminus \{f=0\}$.
Moreover, observe that each point of $\{f=0\}$ is a cluster point of
$W_0$. Suppose by contradiction that $W_0$ extends through a
neighborhood $V$ of $0$ by a complex manifold $W'_0$; then
$\{f=0\}\cap V\subset W'_0$, thus $\{f=0\}\cap V=W'_0\cap V$. This
is a contradiction.
 \eex

\subsection{The unbounded case}

Let $\O\subset\C^n$ be a strictly pseudoconvex domain, and $A\subset
b\O$ an unbounded open subset of $b\O$.

Consider the set
$$\mathcal A\ =\ \left\{A'\Subset b\O\ |\ A'\subset A,\ A'\emph{ domain}\right\}.$$
For an arbitrary $A'\in\mathcal A$ ($bA'=K'$), let $D_{A'}$ be the
compact connected component of $\O\setminus \widehat{K'}$. Set
$$D\ =\ \bigcup_{A'\in\mathcal A} D_{A'}.$$
From Theorem \ref{main} it follows that for every maximally complex
closed ($2m+1$)-dimensional real submanifold $M$ of $A$, there is an
($m+1$)-dimensional complex closed subvariety $W$ of $D$, with
isolated singularities, such that $bW\cap A=M$. So the domain $D$ is
a possible solution of our extension problem.

When $A=b\O$, we may restate the previous result in a more elegant
way. In the same situation as above, consider $$\C^n\subset \C\P^n,\
\C^n=\C\P^n\setminus\C\P^{n-1}_\infty$$ and define the
\emph{principal divisors hull} $\widehat{C}_{\mathcal D}$ of
$C=\oli\O\cap\C\P^{n-1}_\infty$ by
$$\widehat{C}_{\mathcal D}\ =\ \left\{z\in\O\ |\ \forall f\in\mathcal O(\oli\O)\ \oli L_{f,z}\cap C\neq\emptyset \right\},$$
where $\oli L_{f,z}$ is the closure of the
connected component (in $\oli\O$) of the level-set ${\left\{f=f(z)\right\}}$
passing through $z$. Then
$$D=\O\setminus\widehat{C}_{\mathcal D}.$$
Indeed, if $z\in D$, then there are an open subset $A'\subset b\O$
and a function $f\in\mathcal O(\oli\O)$ such that
$\oli L_{f,z}\cap b\O$ is a compact submanifold
of $A'$. In particular $z\not\in\widehat{C}_{\mathcal D}$. Vice
versa, if $z\not\in\widehat{C}_{\mathcal D}$ then there is a
function $g\in\mathcal O(\O')$ ($\O'\supset\O$ domain) such that
$N=\oli L_{g,z}\cap C=\emptyset$, i.e.\ it is a
compact submanifold of $b\O$. By choosing a relatively compact open
subset $A'\subset b\O$ large enough to contain $N$ it follows that
$z\in D_{A'}\subset D$.

\section{Generalization to analytic sets}
Let $\O$, $A$ and $K$ be as before. We want now to consider the extension problem for analytic sets.

Let us recall that if $\mathcal F$ is a coherent sheaf on a domain $U$ in $\C^n$, $x\in U$ and
$$0 \to \mathcal O^{m_k}_x\to \cdots \to\mathcal{O}^{m_0}_x\to\mathcal F_x\to 0$$
is a resolution of $\mathcal F_x$, then the \emph{depth} of $\mathcal F$ at the point $x$ is the integer $p(\mathcal F_x)=n - k$.

We will say that $M\subset A$ is a \emph{$k$-deep trace} of an analytic subset if there are
\begin{itemize}
\item[i)] an open set $U\subset\C^n$ ($U\cap b\Omega=A$);
\item[ii)] an $(m+1)$-dimensional irreducible analytic set $W_M$, whose ideal sheaf $\mathcal I_{W_M}$ has depth at least $k$ at each point of $U$, such that $W_M\cap b\O=M$.
\end{itemize}
In this case, we say that the real dimension of $M$ is $2m+1$.

\bt \label{gener} For any
$(2m+1)$-dimensional $4$-deep trace of analytic subset $M\subset A$, there exists an $(m+1)$-dimensional complex
variety $W$ in $\O \setminus
\widehat K$, such that $bW\cap (A\setminus \widehat K) = M\cap
(A\setminus \widehat K) $.

%In particular, if $K$ is holomorphically convex, i.e.\ $\widehat K=K$, then the extension is in $\O$ and $bW\cap A=M$.
\et

Observe that in this situation we already have a strip $U$ on which the set $M$ extends. So we only need to generalize Lemma \ref{f=0} and the results in Section \ref{sez31}.

\bl \label{f=0gen} Let $z^0\in \O\setminus \widehat K$. Then there
exist an open Stein neighborhood $\O_\alpha \supset \O$ and
$f\in\mathcal O(\O_\alpha)$ such that
\begin{itemize}
  \item [1.]$f(z^0)=0$;
  \item [2.]$\{f=0\}$ is a regular complex hypersurface of $\O_\alpha\setminus\widehat K$;
  \item [3.]$\{f=0\}$ intersects $M$ in a compact set and $W_M$ in an analytic subset (of depth at least $3$).
\end{itemize}
\el

\demo The proof of the first two conditions is exactly the same as before. So, we focus on the third one.

Again, Alexander's Theorem (see \cite[Theorem 3]{A}) implies compactness of the intersection with $M$. Then, we may suppose that $W_M$ is not contained in $\{z_1 =
z_1^0\}$ and, for $\e$ small enough, let $f:\O_\a \to\C$ be the function
$f(z)=h(z) + \e (z_1 - z_1^0)$, where $\O_\a$ and $h$ are as defined in Lemma \ref{f=0}. Consider the stratification of $W_M$ in complex manifolds. By Sard's Lemma, the set of $\e$ for which the intersection of $\{f(z)=0\}$ with a fixed stratum is transversal is open and dense. Hence the set of $\e$ for which the intersection of $\{f(z)=0\}$ with each stratum is transversal is also open and dense, in particular it is non-empty. The conclusion follows. \enddemo

The previous Lemma enables us to extend each analytic subset
$$W_0=W_M\cap\{f=0\}$$
to an analytic set defined on the whole of
$$\O\cap\{f=0\}.$$
Indeed, on a strictly pseudoconvex corona the depth of $W_0$ is at
least $3$ and thus $W_0$ extends in the hole (see e.g.\
\cite{AS,ST}). Obviously the extension lies in $\{f=0\}$.

Observe that, up to a arbitrarily small modification of $b\Omega$ we
can suppose that it intersects each stratum of the stratification of
$W_M$ transversally. In this situation $M$ is a smooth submanifold
with negligible singularities of Hausdorff codimension at least $2$
(see \cite{DH2}).

Again, we consider a generic projection $\pi:\widetilde{U}\to\C^{m}$
and we use holomorphic coordinates $(w',w)$,
$w=(w_1,\ldots,w_{n-m-1})$ on $$\C^n=\C^{m+1}\times\C^{n-m-1}.$$
Keeping the notations used in Section \ref{sez31}, let $V_\tau=
\C^{m+1} \setminus \pi(M_\tau)$.

For $\tau\in U$,
$w'\in\C^{m+1}\setminus\pi(M_\tau)$ and $\alpha\in\N^{n-m-1}$, we
define
$$ I^\alpha(w',\tau)\ \doteqdot\ \int_{(\eta',\eta)\in\ {\rm Reg}(
M_\tau)}\eta^\alpha\omega_{BM}(\eta'-w'), $$
 $\omega_{BM}$ being the
Bochner-Martinelli kernel.

Observe that the previous integral is well-defined and converges. In
fact, $W_\tau=W_M\cap\{f=\tau\}$ is an analytic set and thus, by
Lelong's Theorem, its volume is bounded near the singular locus.
Hence, by Fubini's Theorem, also the regular part of
$M_\tau=W_\tau\cap b\Omega$ has finite volume up to a small
modification of $b\Omega$.

\bl \label{Zaitsevgen} Let $F(w',\tau)$ be the multiple-valued
function which represents $\widetilde M_\tau $ on
$\C^{m}\setminus\pi(M_\tau )$; then, if we denote by
$P^\alpha(F(w',\tau ))$ the sum of the $\alpha^\emph{th}$ powers of the
values of $F(w',\tau )$, the following holds: $$ P^\alpha(F(w',\tau )) =
I^\alpha(w',\tau ). $$ In particular, $F(w',\tau )$ is finite. \el

\demo Let $V_0$ be the unbounded component of $V_\tau $ (where, of course,
$P^\alpha(F(w',\tau )) = 0$). Following \cite
{HL}, it is easy to show that on $V_0$ also $I^\alpha(F(w',\tau )) = 0$. Indeed, if $w'$
is far enough from $\pi({\rm Reg}(M_\tau ))$, then $\beta = \eta^\alpha
\omega_{BM}(\eta' - w')$ is a regular $(m,m-1)$-form on some ball $B_R$ of ${\rm Reg}(M_\tau )$. So, since in $B_R$ there exists $\gamma$
such that $\oli\partial\gamma = \beta$, we may write in the
sense of currents
$$[{\rm Reg}(M_\tau )](\beta) =
[{\rm Reg}(M_\tau )]_{m,m-1}(\oli\partial\gamma) =
\oli\partial[{\rm Reg}(M_\tau )]_{m,m-1}(\gamma).$$

We claim that $\oli\partial[{\rm Reg}(M_\tau )]_{m,m-1}(\gamma)=0$ and, in order to prove this, we first show that $[{\rm Reg}(M_\tau )]$ is a closed current. Indeed, observe that $d[{\rm Reg}(M_\tau )]$ is a flat current, since it is the differential of an ${\rm L^1_{loc}}$ current (see \cite{Fed}). Moreover
$$S={\rm supp}(d[{\rm Reg}(M_\tau )])\subset {\rm Sing}(M_\tau ),$$
hence, denoting by $\dim_{\mathcal H}$ the Hausdorff dimension and by $\mathcal H_s$ the $s$-Hausdorff measure, we have
$$\dim_{\mathcal H}(S)\leq \dim_{\mathcal H} ({\rm sing}(M_\tau ))\leq \dim_{\mathcal H} ({\rm Reg}(M_\tau ))-2$$
and consequently that
$$\mathcal H_{\dim_{\mathcal H}({\rm Reg}(M_\tau ))-1}(S)=0.$$
By Federer's support Theorem (see \cite{Fed}),
this implies that
$$d[{\rm Reg}(M_\tau )]=0.$$

Now, since ${\rm Reg}(M_\tau )$ is maximally complex, $$[{\rm
Reg}(M_\tau )]=[{\rm Reg}(M_\tau )]_{m,m-1} + [{\rm Reg}(M_\tau
)]_{m-1,m}.$$ Since $\oli\partial [{\rm Reg}(M_\tau )]_{m,m-1}$ is
the only component of bidegree ($m,m-2$) of $d[{\rm Reg}(M_\tau )]$
and $d[{\rm Reg}(M_\tau )]=0$ then
$$\oli\partial [{\rm Reg}(M_\tau )]_{m,m-1} = 0.$$

Moreover, since $[{\rm Reg}M_\tau ](\beta)$ is analytic in the
variable $w'$, $[{\rm Reg}M_\tau ](\beta)=0$ for all $w'\in V_0$.

The rest of the proof goes as in Lemma \ref{Zaitsev}.
\enddemo

\bl $P^\alpha(F(w',\tau ))$ is holomorphic in the variable $\tau \in
U\subset\C$, for each $\alpha\in\N^{n-m-1}$. \el \demo The only
difference with the proof for the case of manifolds is the fact that
$I$ is an integration over the regular part of $\Gamma\ast M_\tau $
and not all over $\Gamma\ast M_\tau $. It is easy to see that
Stokes Theorem is valid also in this situation, so the chain of
integrals in Lemma \ref{Lemma36} holds in this case, too.
\enddemo

The rest of the proof of Theorem \ref{gener} goes as in the proof of Theorem \ref{main} (see Section \ref{sez31}).
%\end {}

\end{document}